\documentclass{amsart}
\usepackage{amssymb}
\usepackage{mathrsfs}
\usepackage{graphicx}
\usepackage{comment}

\def\S{\mathcal S}

\newcommand{\nni}{\noindent}
\newcommand{\be}{\begin{equation}}
\newcommand{\ee}{\end{equation}}
\newcommand{\ba}{\begin{align}}
\newcommand{\ea}{\end{align}}

\newcommand{\abs}[1]{\lvert#1\rvert}

\newtheorem{example}{Example}[section]
\newtheorem{theorem}{Theorem}[section]

{\begin{list}{}{%
\settowidth{\labelwidth}{\textsf{{\it #1.}}}%
\setlength{\labelsep}{4mm}%
\setlength{\leftmargin}{\labelwidth}%
\addtolength{\leftmargin}{\labelsep}%
}}%
{\end{list}}

\def\beq{\begin{equation}}\def\enq{\end{equation}}

{\begin{list}{}{%
\settowidth{\labelwidth}{\textsf{{\it #1.}}}%
\setlength{\labelsep}{2mm}%
\setlength{\leftmargin}{\labelwidth}%
\addtolength{\leftmargin}{\labelsep}%
\addtolength{\leftmargin}{4mm}%
\setlength{\itemsep}{6pt}%
\setlength{\listparindent}{0pt}%
\setlength{\topsep}{3pt}%
}}%

\title[Integer group determinants for SmallGroup(16,8)]{The integer group determinants for the semidihedral group of order 16}

\author[H. Bautista Serrano]{Humberto Bautista Serrano}

\author[B. Paudel]{Bishnu Paudel}
\address{ Department of Mathematics\\
         Kansas State University\\
         Manhattan, KS 66506, USA}
\email{humbertb@ksu.edu, bpaudel@ksu.edu, pinner@math.ksu.edu}

\author[C. Pinner]{Chris Pinner}


\date{\today}
\begin{document}

\begin{abstract}
We obtain a complete description of the integer group determinants for SmallGroup(16,8), 
the semidihedral group of order 16. 

While this paper was in preparation, a complete descriptions  for this group was 
independently obtained  by Yuka Yamaguchi and Naoya Yamaguchi in \cite{Yamaguchi9} (and the other remaining 
group of order 16 in \cite{Yamaguchi8}). We offer our version here  for anyone interested in a slightly different  approach.

\end{abstract}

\maketitle

\section{Introduction} 

It is a non-trivial problem to determine the integer values are taken by the group 
determinant when the entries are all integers, a question posed by Olga Taussky-Todd \cite{TausskyTodd} at the meeting of the American Mathematical Society in Hayward, California, in April 1977. A complete description has been obtained
for all groups of order at most 15 (see \cite{Newman1,Laquer,dihedral,smallgps,bishnu1}
and for 12 of the 14 groups of order 16 (the 5 abelian groups in \cite{Yamaguchi1,Yamaguchi2,Yamaguchi3,Yamaguchi4,Yamaguchi5}), and for 7 of the  9 non-abelian groups in \cite{dihedral,ZnxH,Q16,Yamaguchi6,Yamaguchi7, Humb}).
We recall the definition of the group determinant, noting that we shall think of the integer group
determinant as being defined on the elements in the group ring $\mathbb Z[G]$
$$ \mathcal{D}_G\left( \sum_{g\in G} a_g g \right):=\det\left( a_{gh^{-1}}\right) .$$

In this paper we determine  the integer group determinants for the 
semidihedral or quasidihedral group of order 16
$$SD_{16} \text{ or } QD_{16} =   \langle X,Y\; | \; X^8=Y^2=1,\; YXY=X^3\rangle. $$
This is SmallGroup(16,8), using  the group identification from
GAP's small group library.
Notice that we have order 8 dihedral and quaternion subgroups $D_8=\langle X^2, Y\rangle$ and
$Q_8=\langle X^2,XY\rangle$ respectively.

\begin{theorem} The even integer group determinants for $SD_{16}$ are exactly the multiples of $2^{10}$.

The odd integer determinants are the integers $n\equiv 1$ mod 8 plus those $n\equiv 5$ mod 8 of the form $n=mp^2$ with $m\equiv 5 $ mod 8 and $p$ a prime which is 3 mod 8.

\end{theorem}
Recently Yuka and Naoya Yamaguchi have independently obtained this result \cite{Yamaguchi9}, plus the determinants for the remaining group of order 16 in \cite{Yamaguchi8}. The relationship between the integer group determinants for 
the 14 groups of order 16 is nicely illustrated in \cite{Yamaguchi9}.

\section{A formula for the semidihedral group determinant}

Frobenius \cite{Frob} observed that the group determinant can be factored 

$$ \mathcal{D}_G\left( \sum_{g\in G} a_g g \right)=\prod_{\rho\in \hat{G}} \det \left(\sum_{g\in G} a_g \rho(g) \right)^{\deg \rho}, $$
where $\hat{G}$ denotes a complete set of non-isomorphic irreducible group  representations  for the group (see for example \cite{Conrad} or \cite{book}).

We use this to give the formula for the general semidihedral (also called quasidihedral) group of order $2^n$, $n\geq 4$
$$ SD_{2^n}=\left\langle X,Y\; |\; X^{2^{n-1}}=Y^2=1,\; YXY=X^{2^{n-2}-1}\right\rangle. $$

From the group presentation, an element in the group ring $\mathbb Z[SD_{2^n}]$ can be written in the  form
\be \label{Fform} F(X,Y)= f(X) + Y g(X) \ee
where $f$ and $g$ are integer polynomials of degree $2^{n-1}-1$.

From the group presentation we know that the group characters must satisfy
$\chi(Y)^2=\chi(X)^{2^{n-1}}=\chi(X)^{2^{n-2}-2}=1$, that is, we have four
characters:
$$ \chi(X),\chi(Y)=\pm 1. $$
These contribute
\be \label{defM} M=F(1,1)F(1,-1)F(-1,1)F(-1,-1)=\left( f(1)^2-g(1)^2 \right)(f(-1)^2-g(-1)^2)\ee
to the group determinant. This is the $\mathbb Z_2\times \mathbb Z_2$ determinant of $F(X,Y)$.

This leaves $2^{n-2}-1$ degree two representations to find. Now if $\lambda^{2^{n-1}}=1$, then $\lambda^{(2^{n-2}-1)^2}=\lambda$ and
$$ \rho(X)=\begin{pmatrix} \lambda & 0 \\ 0 & \lambda^{2^{n-2}-1} \end{pmatrix},\quad \rho(Y)=\begin{pmatrix} 0 & 1 \\ 1 & 0 \end{pmatrix} $$
satisfies the group relations, with
\begin{align*} \det \left(\rho(F(X,Y)\right) & = \det \begin{pmatrix}  f(\lambda) & g(\lambda^{2^{n-2}-1}) \\ g(\lambda) & f(\lambda^{2^{n-2}-1})\end{pmatrix}\\
& =f(\lambda)f(\lambda^{2^{n-2}-1})-g(\lambda)g(\lambda^{2^{n-2}-1}):=k(\lambda). 
\end{align*}
The complex $\lambda$ give us  $2^{n-1}-2$ additional quadratic factors and:
$$\mathcal{D}_{SD_{2^n}}(F(X,Y))= \prod_{\lambda^{2^{n-1}}=1}f(\lambda)f(\lambda^{2^{n-2}-1})-g(\lambda)g(\lambda^{2^{n-2}-1}). $$
That is, the $SD_{2^n}$ determinant of $f(x)+yg(x)$ is the $\mathbb Z_{2^{n-1}}$ determinant of 
$$f(x)f(x^{2^{n-2}-1})-g(x)g(x^{2^{n-2}-1}). $$
Notice that the complex roots $\lambda$ have $k(\lambda)=k(\lambda^{2^{n-2}-1})$ so we will repeat these  factors twice (that is we needed only half the complex $\rho$, with the resulting factor squared). For the $\lambda^{2^{n-2}}=1$ we plainly have 
$$k(\lambda)=f(\lambda)f(\lambda^{-1})-g(\lambda) g(\lambda^{-1}), $$
and when $\lambda$ is a primitive $2^{n-1}$st root of unity
$$k(\lambda)=f(\lambda)f(-\lambda^{-1})-g(\lambda) g(-\lambda^{-1}). $$
Hence we can write
$$\mathcal{D}_{SD_{2^n}}(F(X,Y))= M \prod_{j=2}^{n-1}A_j^2 $$
with $M$ as in \eqref{defM}, for $2\leq j\leq n-2$ the integers
$$ A_j=\prod_{\stackrel{i=1}{i \text{ odd }}}^{2^{j-1}} f(\omega_{2^j}^i)f(\omega_{2^j}^{-i})- g(\omega_{2^j}^i)g(\omega_{2^j}^{-i}),$$
where $\omega_t:=e^{2\pi i/t}$ denotes a primitive $t$th root of unity, and integer
$$ A_{n-1} =\prod_{\stackrel{i=1}{i \text{ odd}}}^{2^{n-3}} \abs{f(\omega_{2^{n-1}}^i)f(-\omega_{2^{n-1}}^{-i})- g(\omega_{2^{n-1}}^i)g(-\omega_{2^{n-1}}^{-i})}^2. $$
Notice that $M\prod_{j=2}^{n-2}A_j^2$ is the dihedral,  $D_{2^{n-1}}$,
determinant for $F(X,Y)$ in \cite{dihedral}.

\subsection{ Restricting to the case of $SD_{16}$}

When $n=4$ this  becomes
\be \label{deffF} F(X,Y)=f(X)+Yg(X),\quad f(x)=\sum_{i=0}^7 a_ix^i,\;g(x)=\sum_{i=0}^7 b_ix^i\in \mathbb Z[x], \ee
and
\begin{align*} \mathcal{D}_{SD_{16}}(F(X,Y)) & =\prod_{j=0}^7 f(\omega_8^j)f(\omega_8^{3j})- g(\omega_8^j)g(\omega_8^{3j})\\ 
& = MA_2^2A_3^2 
\end{align*}
with $M$ as in \eqref{defM} and integers
$$ A_2= f(i)f(-i)-g(i)g(-i) $$
and
$$ A_3=\left( f(\omega_8)f(\omega_8^3)-g(\omega_8)g(\omega_8^3)\right)
\left( f(\omega_8^{5})f(\omega_8^7)-g(\omega_8^5)g(\omega^7)\right),$$
where $\omega_8=\frac{\sqrt{2}}{2}(1+i)$, $\omega_8^3=-\overline{\omega_8}$, $\omega_8^5=-\omega_8$, $\omega^7=\overline{\omega}_8$.

\subsection{ The groups $M_{2^n}$}

The representations are very similar for the family
$$ M_{2^{2^n}}=\left\langle X,Y\; |\; X^{2^{n-1}}=Y^2=1,\; YXY=X^{2^{n-2}+1}\right\rangle. $$
Again we can  take $F(X,Y)$ of the form \eqref{Fform}. This time
we have more characters; $\chi(Y)=\pm 1$ and $\chi(X)$ a $2^{n-2}$th root of unity,
contributing
\be \label{defM1} M_1:=\prod_{x^{2^{n-2}}=1, \; y^2=\pm 1} F(x,y)
= \prod_{x^{2^{n-2}}=1} f(x)^2-g(x)^2, \ee
to the group determinant; that is the $\mathbb Z_{2^{n-2}}\times \mathbb Z_2$ determinant of $F(X,Y)$.

The degree two representations are as before except we take 
$$\rho(X)=\begin{pmatrix} \lambda & 0 \\ 0 & \lambda^{2^{n-2}+1} \end{pmatrix}, $$
with the new factors coming from the primitive $2^{n-1}$st roots of unity. Hence
$$ \mathcal{D}_{M_{2^n}}=\prod_{x^{2^{n-1}}=1} f(x)f(x^{2^{n-2}+1})- g(x)g(x^{2^{n-2}+1}) =M_1A^2, $$
with $M_1$ as in \eqref{defM1} and
$$ A=\prod_{\stackrel{1\leq j< 2^{n-2}}{j \text{ odd }}} f\left(\omega_{2^{n-1}}^j\right)f\left(-\omega_{2^{n-1}}^j\right)- g\left(\omega_{2^{n-1}}^j\right)g\left(-\omega_{2^{n-1}}^j\right). $$

For example, for $n=4$ and $M_{16}=$SmallGroup(16,6), we take  $F(X,Y)$ of the form \eqref{deffF} and get
$$ \mathcal{D}_{M_{16}}(F(X,Y))= M_1A^2$$
with integers
$$ M_1=\prod_{x=\pm 1,\pm i, y=\pm 1} F(x,y) $$
and
$$ A=\left( f(\omega_8)f(-\omega_8) - g(\omega_8)g(-\omega_8) \right)\left( f(\omega_8^3)f(-\omega_8^3) - g(\omega_8^3)g(-\omega_8^3) \right).
$$

\section{Achieved values for $SD_{16}$}

\noindent
{\bf Even Values.} We achieve all multiples of $2^{10}$.

Writing
$$ h(x):=(x+1)(x^2+1)(x^4+1), $$
we get $2^{12}m$ from
$$ F(x,y)=(1+x+x^2+x^3-x^6-x^7-mh)+y(1+x-x^3-x^4-x^5-x^7-mh), $$
we get $2^{11}(2m+1)$ from
$$ (1+x+x^2+x^3+x^4+x^5+x^6-x^7+mh)+y(1+x-x^5+x^7+mh), $$
with $2^{10}(-1+4m)$ from
$$ \left((1+x)(1+x^2)-mh\right) + y \left((1+x+x^2)(1-x^5)-mh\right),  $$
and $2^{10}(1+4m)$ from
$$  (1+x+x^2+x^3+x^4-x^7+mh)+y(1+x-x^3-x^4-x^5+x^7+mh).  $$

\vskip0.1in
\noindent
{\bf Odd Values $n\equiv 1 $ mod 8.} 

We get $1+16m$ from $F(x,y)=(1+mh)+ymh$ and $(16m-7)$ from
$$F(x,y)=(1+x+x^2+x^3+x^4-x^6-mh)+y(1+x+x^5-mh). $$

\vskip0.1in
\noindent
{\bf Odd Values $mp^2\equiv 5 $ mod 8.} 

Suppose that $p\equiv 3$ mod 8. Then $\left(\frac{-2}{p}\right)=1$ and $p$ splits in $\mathbb Z[\sqrt{2}i]$. Hence there are integers $U$, $V$ with $U^2+2V^2=p,$
with $U,V$ necessarily odd. 

Changing the sign of $U$ or $V$ as needed, we write  $U=4k+1$, $V=4s+1$ 
and get $(5+16m)p^2$ from
$$F(x,y)=(x+x^2-(k-sx^2)(1-x^4)+mh)+y(1+x+x^2+(s+kx^2)(1-x^4)+mh), $$
writing $U=4s+1$ and $V=4k-1$ we get $(16m-3)p^2$ from
$$  F(x,y)= (1+x+(s+kx^2  )(1-x^4)-mh)+y(x+(k-sx^2)(1-x^4)-mh).$$

\section{Values must be of the stated form}

It was shown in  \cite{dihedral} that an odd  $D_8$ determinant must be $1$ mod 4 and an even determinant  be divisible by $2^8$. 
Since $\mathcal{D}_{SD_{16}}(F)=MA_2^2A_3^2$ where $MA_2^2$ is the $D_8$ determinant of $F$, 
we see that an odd $SD_{16}$ determinant must be 1 mod 4. Since
$MA_2^2\equiv F(1,1)^8$ mod 2 and $A_3\equiv F(1,1)^4$ mod 2, we see that an even determinants must have $F(1,1)$ even, $2^8\mid MA_2^2$, $2^2\mid A_3^2$ and $2^{10}\mid \mathcal{D}_{SD_{16}}(F)$; all of which we can achieve.

Since we achieve the 1 mod 8, it only remains to show that the 5 mod 8 are of the stated form. 
Suppose that $\mathcal{D}_{SD_{16}}(F)\equiv 5$ mod 8. As it is odd we know $f(1)$ and $g(1)$ have opposite parity. We suppose that $f(1),f(-1)$ are odd and $g(1)$, $g(-1)$
are even (else switch the roles of $f$ and $g$). Since $(A_2A_3)^2\equiv 1$ mod 8 we must
have one of $g(1)^2,g(-1)^2\equiv 4$ mod 8 and the other 0 mod 8. That is,
\be \label{diff}  \frac{1}{4}(g(1)^2-g(-1)^2)\equiv 1 \text{ mod } 2. \ee
We write
\begin{align*}
\alpha_0=a_0-a_4, \quad \alpha_1=a_1-a_5,\quad & \alpha_2=a_2-a_6, \quad \alpha_3=a_3-a_7,\\
\gamma_0=a_0+a_4, \quad \gamma_1=a_1+a_5,\quad & \gamma_2=a_2+a_6, \quad \gamma_3=a_3+a_7,
\end{align*}
and $\beta_0=b_0-b_4,\ldots,\beta_3=b_3-b_7$ and $\delta_0=b_0+b_4,\ldots ,\delta_3=b_3+b_7.$
Then 
\begin{align*} f(\omega_8) & = \alpha_0 +\alpha_1 \omega_8 + \alpha_2 i -\alpha_3\overline{\omega}_8\\
 & = \left( \alpha_0 +\frac{\sqrt{2}}{2}i(\alpha_1+\alpha_3)\right) +\left(i\alpha_2+\frac{\sqrt{2}}{2}(\alpha_1-\alpha_3)\right), \\
f(-\overline{\omega}_8) & = \left( \alpha_0 +\frac{\sqrt{2}}{2}i(\alpha_1+\alpha_3)\right) -\left(i\alpha_2+\frac{\sqrt{2}}{2}(\alpha_1-\alpha_3)\right), 
\end{align*}
and
$$ f(\omega_8)f(-\overline{\omega}_8)= U_1+ \sqrt{2}i  V _1$$
with
\begin{align*}
U_1 &= \alpha_0^2-\alpha_1^2+\alpha_2^2-\alpha_3^2,\\
V_1& = \alpha_0\alpha_1+\alpha_0\alpha_3-\alpha_1\alpha_2+\alpha_2\alpha_3.
\end{align*}
Plainly
$$ U_1\equiv \gamma_0^2+\gamma_1^2+\gamma_2^2+\gamma_3^2 \equiv f(1)^2\equiv 1 \text{ mod } 2. $$
Likewise 
$$ V_1 \equiv \gamma_0\gamma_1+\gamma_0\gamma_3+\gamma_1\gamma_2+\gamma_2\gamma_3 \text{ mod }2. $$
Since $f(1)=\gamma_0+\gamma_1+\gamma_2+\gamma_3$ is odd, we must have one of 
$\gamma_0,\ldots ,\gamma_3$ odd and three even, or three odd and one even. In either case we get that $V_1$ is even. 

Likewise we can write
$$ g(\omega_8)g(-\overline{\omega}_8)= U_2 + i\sqrt{2} V_2$$
with
\begin{align*}
U_2 & \equiv g(1)^2 \equiv 0 \text{ mod } 2,\\
V_2 & \equiv  \delta_0\delta_1+\delta_0\delta_3+\delta_1\delta_2+\delta_2\delta_3 \text{ mod }2. 
\end{align*}
But the latter expression is \eqref{diff} and hence $V_2$ is odd. Thus $(U_1-U_2)$ 
and $(V_1-V_2)$ are odd and 
$$ A_3 = (U_1-U_2)^2+2(V_1-V_2)^2 \equiv 3 \text{ mod } 8. $$
Notice that 
$$\left(\frac{-2}{p}\right) =\begin{cases} 1 & \text{if $p\equiv 1$ or 3 mod 8,}\\
 -1 & \text{ if $p\equiv 5$ or 7 mod 8,} \end{cases} $$
so that primes $5$ or $7$ mod 8 in $A_3$ must be squared. That is, $A_3$ must contain a prime
$p\equiv 3$ mod 8 and $\mathcal{D}_{SD_{16}}(F)=mp^2$ with $m\equiv 5$ mod 8 
and $p\equiv 3$ mod 8. \qed

\end{document}